\documentclass[11pt]{article} \makeatother
\usepackage{amsfonts, amsmath, amssymb}
\usepackage{graphicx}
\usepackage{subfigure}
\usepackage{mathrsfs}
\allowdisplaybreaks[4]
\newtheorem{theo}{Theorem}[section]
\newtheorem{lemma}[theo]{Lemma}

\newcommand{\be}{\begin{equation}}
\newcommand{\ee}{\end{equation}}
\newcommand\bes{\begin{eqnarray}} \newcommand\ees{\end{eqnarray}}
\newcommand{\bess}{\begin{eqnarray*}}
\newcommand{\eess}{\end{eqnarray*}}
\numberwithin{equation}{section}

\AtBeginDocument{{\noindent\small
}
}

\begin{document}

\date{}
\setlength{\baselineskip}{17pt}{\setlength\arraycolsep{2pt}
%

\begin{center} {\bf\Large Well-posedness and general decay of solution }\\[2mm]
 {\bf\Large for a transmission problem with viscoelastic term and delay}\\[4mm]
  {\large Danhua Wang\footnote{E-mail: matdhwang@yeah.net}, Gang Li and Biqing Zhu}\\[1mm]
{\small  College of Mathematics and Statistics, Nanjing University of
Information Science and Technology, Nanjing 210044, China} \\[2mm]
\end{center}

\begin{abstract}
In this paper, we consider a transmission problem in a bounded domain with a viscoelastic term and a delay term.
Under appropriate hypothesis on the relaxation function and the relationship between the weight of the damping and the weight of the delay, we prove the well-posedness result by using Faedo-Galerkin method. By introducing suitable Lyaponov functionals, we establish a general decay result, from which the exponential and polynomial types of decay are only special cases.
\end{abstract}



\section{Introduction}\label{s1}
\setcounter{equation}{0}

In this paper,  we study the transmission system with
a viscoelastic term and a delay term
 \bes\left\{\begin{array}{ll}
\displaystyle u_{tt}(x,t)-au_{xx}(x,t)+\int_{0}^{t}g(t-s)u_{xx}(x,s){\rm d}s\\\medskip\displaystyle
\quad\quad\quad\quad\quad\quad\quad\quad\quad  +\mu_{1}u_{t}(x,t)+\mu_{2}u_{t}(x,t-\tau)=0,
&(x,t)\in\Omega\times(0,+\infty),
\medskip\\
  \displaystyle v_{tt}(x,t)-bv_{xx}(x,t)=0,& (x,t)\in(L_{1},L_{2})\times(0,+\infty),
 \end{array}\right.\label{1.1}
 \ees
under the boundary and transmission conditions
 \bes\left\{\begin{array}{ll}
\displaystyle u(0,t)=u(L_{3},t)=0,
\medskip\\\medskip
\displaystyle u(L_{i},t)=v(L_{i},t), &i=1,2,
\medskip\\
  \displaystyle \left(a-\int_{0}^{t}g(s){\rm d}s\right)u_{x}(L_{i},t)=bv_{x}(L_{i},t), &i=1,2,
 \end{array}\right.\label{1.2}
 \ees
 and the initial conditions
  \bes\left\{\begin{array}{ll}
\displaystyle u(x,0)=u_{0}(x),\quad u_{t}(x,0)=u_{1}(x), &x\in\Omega,
\medskip\\\medskip
\displaystyle u_{t}(x,t-\tau)=f_{0}(x,t-\tau), &x\in\Omega,\quad t\in[0,\tau],
\medskip\\
  \displaystyle v(x,0)=v_{0}(x),\quad v_{t}(x,0)=v_{1}(x), &x\in(L_{1},L_{2}),
 \end{array}\right.\label{1.3}
 \ees
 where $0<L_{1}<L_{2}<L_{3}$, $\Omega=(0,L_{1})\cup(L_{2},L_{3})$, $a$, $b$, $\mu_{1}$, $\mu_{2}$ are positive constants, and
 $\tau>0$ is the delay.

\begin{figure}[!ht]\label{1.5}
\begin{center}
\includegraphics[width=0.6\textwidth]{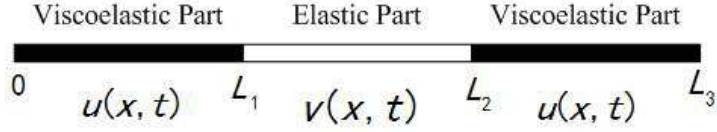}
\caption{The configuration.}\label{figure1}
\end{center}
\end{figure}
 The problems like \eqref{1.1}-\eqref{1.3} related to the wave propagation over a body are called transmission problems, which consists of two different types of materials: the elastic part and the viscoelastic part.

In recent years, many authors have investigated wave equations with viscoelastic damping and showed that the dissipation produced by the viscoelastic part can produce the decay of the solution, see  \cite{bm2004,bm2006,c2001,cds2002,h2011,l2014,lc2014,ls2014,m2008,tp2013,w2013,w20132} and the references therein. For example,
  Cavalcanti et al. \cite{cds2002} studied the following equation:
\begin{equation*}
u_{tt}-\Delta u+\int_{0}^{t}g(t-\tau)\Delta u(\tau){\rm d}\tau+a(x)u_{t}+|u|^{\gamma}u=0,\quad {\rm in}\quad \Omega\times(0,\infty),
\end{equation*}
where $a: \Omega\rightarrow \mathbb{R}_{+}$. Under the conditions
that $a(x)\geq a_{0}>0$ on $\omega\subset\Omega$, with $\omega$ satisfying some geometry restrictions and
\begin{equation*}
-\xi_{1}g(t)\leq g'(t)\leq-\xi_{2}g(t),\quad t\geq 0,
\end{equation*}
the authors showed the exponential decay. Then Berrimi and Messaoudi \cite{bm2004} proved the same
result under weaker conditions on both $a$ and $g$.
 Berrimi and Messaoudi \cite{bm2006}  considered the equation
\begin{equation*}
u_{tt}-\Delta u+\int_{0}^{t}g(t-\tau)\Delta u(\tau){\rm d}\tau=|u|^{\gamma}u,\quad {\rm in}\quad \Omega\times(0,\infty),
\end{equation*}
with only the viscoelastic dissipation and proved
that the solution energy  decays
exponentially or polynomially depending on the rate of the decay of the relaxation function $g$.
In all previous works, the rates of decay of relaxation functions were either exponential or polynomial type.
For a wider class of relaxation functions,
Messaoudi \cite{m2008} investigated the following viscoelastic equation:
\begin{equation*}
u_{tt}-\Delta u+\int_{0}^{t}g(t-\tau)\Delta u(\tau){\rm d}\tau=0,\quad {\rm in}\quad \Omega\times(0,\infty),
\end{equation*}
in a bounded domain, and established a more
general decay result, from which the usual exponential and polynomial decay rates
are only special cases.
Afterwards, Han and Wang \cite{h2011} studied the nonlinear viscoelastic equation
\begin{equation*}
u_{tt}-\Delta u+\int_{0}^{t}g(t-\tau)\Delta u(\tau){\rm d}\tau+|u|^{k}\partial j(u_{t})=|u|^{p-1}u,\quad {\rm in}\quad \Omega\times(0,T).
\end{equation*}
They obtained the global existence of generalized solutions, weak solutions for
the equation. In addition, the finite time blow-up of weak solutions is established provided that the initial energy is negative
and the exponent $p$ is greater than the critical value.

It is well known that delay effects, which arise in many practical problems, may be sources of instability.
Hence, the control of PDEs with time delay effects has become an active area of research in recent years.
For example, it was proved in
\cite{dl1986,np2006} that an arbitrarily small delay may destabilize
a system which is uniformly asymptotically stable in the absence of delay unless
additional conditions or control terms were used.
A boundary stabilization problem for the wave equation with interior delay studied in \cite{anp2010}. The authors proved an exponential stability result under some Lions geometric condition.
Kirane and Said-Houari \cite{ks2011} considered the viscoelastic wave equation with a delay
$$
u_{tt}(x,t)-\Delta u(x,t)+\int_{0}^{t}g(t-s)\Delta u(x,t-s){\rm d}s+\mu_{1}u_{t}(x,t)+\mu_{2}u_{t}(x,t-\tau)=0,\quad {\rm in}\quad \Omega\times(0,\infty),
$$
where $\mu_{1}$ and $\mu_{2}$ are positive constants. They established a general energy decay result under the condition
that $0\le \mu_{2}\le \mu_{1}$.
Later, Liu \cite{l2013} improved this result by considering the equation with a time-varying delay term, with not necessarily positive coefficient $\mu_{2}$ of the delay term.

Transmission  problems
related to \eqref{1.1}-\eqref{1.3} have also been extensively studied. Bastos and Raposo \cite{br2007}
investigated the transmission problem with frictional damping
and  showed the well-posedness and exponential stability of the total energy.
 Mu\~noz Rivera and
Portillo Oquendo \cite{mp2000} considered the transmission problem of viscoelastic waves
and proved that    the dissipation produced by the viscoelastic part can produce exponential decay
of the solution, no matter how small its size is.
Bae \cite{b2010} studied the transmission problem,
in which one component is clamped and the other is in a viscoelastic
fluid producing a dissipative mechanism on the boundary, and established a decay result
which depends on the rate of the decay of the relaxation function.


Motivated by the above results, we intend to consider the well-posedness and the general decay result of
problem \eqref{1.1}-\eqref{1.3} under some hypothesis in this paper.
The main difficulty we encounter here arises from
the simultaneous appearance of the viscoelastic term and the delay term.
Our first intention is to study the well-posedness
of problem \eqref{1.1}-\eqref{1.3} by making use of Faedo-Galerkin procedure, that is Faedo-Galerkin approximation together with energy estimates.
For asymptotic behavior, we prove a general decay result from which the exponential and polynomial types of decay are only special cases by introducing suitable Lyaponov functionals.

The paper is organized as follows. In Section \ref{s2}, we give some materials needed for our work and state our main results. In Section \ref{s3}, we prove the well-posedness of the problem. The general decay result is proved in Section \ref{s4}.

\section{Preliminaries and main results}\label{s2}
\setcounter{equation}{0}

In this section, we present some materials that shall be used in order to prove our main results.
Let us first introduce the following notations:
\begin{align*}
&(g*h)(t):=\int_{0}^{t}g(t-s)h(s){\rm d}s,\\
&(g\diamond h)(t):=\int_{0}^{t}g(t-s)|h(t)-h(s)|{\rm d}s,\\
&(g\square h)(t):=\int_{0}^{t}g(t-s)|h(t)-h(s)|^2{\rm d}s.
\end{align*}
Note that the sign of $(g\square h)(t)$ depends solely on the sign of $g$. We easily see that the above operators satisfy
\begin{align*}
&(g*h)(t)=\left(\int_{0}^{t}g(s){\rm d}s\right)h(t)-(g\diamond h)(t),\\
&|(g\diamond h)(t)|^2\leq\left(\int_{0}^{t}|g(s)|{\rm d}s\right)(|g|\square h)(t).
\end{align*}

\begin{lemma}\label{le2.1}
For any $g,h\in C^{1}(\mathbb{R})$, the following equation holds
\begin{equation*}
2[g*h]h'=g'\square h-g(t)|h|^2-\frac{d}{dt}\left\{g\square h-\left(\int_{0}^{t}g(s){\rm d}s\right)|h|^2\right\}.
\end{equation*}
\end{lemma}

For the relaxation function $g$, we assume

(G1) $g$: $\mathbb{R}_{+}\rightarrow\mathbb{R}_{+}$ is a $C^{1}$ function satisfying
\begin{equation*}
g(0)>0,\quad 0<\beta(t):=a-\int_{0}^{t}g(s){\rm d}s\quad{\rm and}\quad 0<\beta_{0}:=a-\int_{0}^{\infty}g(s){\rm d}s.
\end{equation*}

(G2) There exists a nonincreasing differentiable function $\xi(t)$: $\mathbb{R}_{+}\rightarrow\mathbb{R}_{+}$ such that
\begin{equation*}
g'(t)\leq -\xi(t) g(t), \quad \forall t\geq 0 \quad{\rm and}\quad \int_{0}^{\infty}\xi(t){\rm d}t=+\infty.
\end{equation*}
These hypotheses imply that
\begin{equation}
\beta_{0}\leq\beta(t)\leq a.
\end{equation}

As in \cite{np2006}, we introduce the following variable:
\begin{equation*}
z(x,\rho,t)=u_{t}(x,t-\tau\rho),\quad (x,\rho,t)\in \Omega\times(0,1)\times(0,\infty).
\end{equation*}
Then the above variable $z$ satisfies
\begin{equation}
\tau z_{t}(x,\rho,t)+z_{\rho}(x,\rho,t)=0,\quad (x,\rho,t)\in \Omega\times(0,1)\times(0,\infty).
\end{equation}
Thus, system \eqref{1.1} becomes
\bes\left\{\begin{array}{ll}
\displaystyle u_{tt}(x,t)-au_{xx}(x,t)+g*u_{xx}\\\medskip\displaystyle
\quad\quad\quad\quad\quad\quad\quad\quad\quad  +\mu_{1}u_{t}(x,t)+\mu_{2}z(x,1,t)=0,&(x,t)\in\Omega\times(0,+\infty),
\\
  \displaystyle v_{tt}(x,t)-bv_{xx}(x,t)=0,& (x,t)\in(L_{1},L_{2})\times(0,+\infty),
\medskip\\
\displaystyle  \tau z_{t}(x,\rho,t)+z_{\rho}(x,\rho,t)=0,& (x,\rho,t)\in\Omega\times(0,1)\times(0,+\infty),
 \end{array}\right.\label{2.5}
 \ees
 and the boundary and transmission conditions \eqref{1.2}  becomes
\bes\left\{\begin{array}{ll}
\displaystyle u(0,t)=u(L_{3},t)=0,
\medskip\\\medskip
\displaystyle u(L_{i},t)=v(L_{i},t),\quad\quad\quad\quad\quad\quad\quad\quad\quad\quad i=1,2,\quad &t\in(0,+\infty),
\medskip\\\medskip
  \displaystyle \left(a-\int_{0}^{t}g(s){\rm d}s\right)u_{x}(L_{i},t)=bv_{x}(L_{i},t),\quad i=1,2,\quad &t\in(0,+\infty),
 \medskip\\\medskip
  \displaystyle z(x,0,t)=u_{t}(x,t), &(x,t)\in\Omega\times(0,+\infty),
 \medskip\\
  \displaystyle z(x,1,t)=f_{0}(x,t-\tau), &(x,t)\in\Omega\times(0,\tau).
 \end{array}\right.\label{2.6}
 \ees

Similar to \cite{r2008}, we denote the Hilbert spaces
\begin{align*}
\mathcal{V}=\bigg\{&(u,v)\in H^{1}(\Omega)\cap H^{1}(L_{1},L_{2}):u(0,t)=u(L_{3},0)=0,u(L_{i},t)=v(L_{i},t),\bigg.\\
&\bigg.\left(a-\int_{0}^{t}g(s){\rm d}s\right)u_{x}(L_{i},t)=bv_{x}(L_{i},t), i=1,2\bigg\}
\end{align*}
and
\begin{equation*}
\mathcal{L}^2=L^{2}(\Omega)\times L^{2}(L_{1},L_{2}).
\end{equation*}
Then the existence result reads as follows:
\begin{theo}\label{th2.1}
Assume that $\mu_{2}\leq\mu_{1}$, $(G1)$ and $(G2)$ hold. Then given $(u_{0},v_{0})\in \mathcal{V}$, $(u_{1},v_{1})\in \mathcal{L}^2$, and $f_{0}\in L^2((0,1),H^{1}(\Omega))$, there exists a unique weak solution $(u,v,z)$ of problem \eqref{2.5}-\eqref{2.6} such
that
\begin{align*}
&(u,v)\in C(0,\infty;\mathcal{V})\cap C^{1}(0,\infty;\mathcal{L}^2),\\
&z\in C(0,\infty;L^2((0,1),H^{1}(\Omega))).
\end{align*}
\end{theo}

For any regular solution
of \eqref{1.1}-\eqref{1.3}, we define the energy as
\begin{align}\label{3.1}
E(t)=&\frac{1}{2}\int_{\Omega}u_{t}^2(x,t){\rm d}x+\frac{1}{2}\beta(t)\int_{\Omega}u_{x}^2(x,t){\rm d}x+\frac{1}{2}\int_{\Omega}(g\square u_{x}){\rm d}x\nonumber\\
&+\frac{1}{2}\int_{L_{1}}^{L_{2}}\left[v_{t}^2(x,t)+bv_{x}^2(x,t)\right]{\rm d}x+\frac{\zeta}{2}\int_{\Omega}\int_{0}^{1}z^2(x,\rho,t){\rm d}\rho{\rm d}x,
\end{align}
where $\zeta$ is a positive constant such that
\begin{equation}\label{2.9}
\tau\mu_{2}<\zeta<\tau(2\mu_{1}-\mu_{2}).
\end{equation}
Our decay result reads as follows:
\begin{theo}\label{th3.1}
Let $(u,v)$ be the solution of problem \eqref{1.1}-\eqref{1.3}. Assume that $\mu_{2}<\mu_{1}$, $(G1)$, $(G2)$ and
\begin{equation}\label{3.30}
b>\frac{4(L_{2}-L_{1})}{L_{1}+L_{3}-L_{2}}\beta_{0},\quad a>\frac{4(L_{2}-L_{1})}{L_{1}+L_{3}-L_{2}}\beta_{0}
\end{equation}
hold,
then there exists constants $\gamma_{0},\gamma_{2}>0$ such that, for all $t\in\mathbb{R}_{+}$ and for all $\gamma_{1}\in(0,\gamma_{0})$,
\begin{equation}\label{3.26}
E(t)\leq \gamma_{2}e^{-\gamma_{1}\int_{0}^{t}\xi(s){\rm d}s}.
\end{equation}
\end{theo}

\section{Well-posedness of the problem}\label{s3}
\setcounter{equation}{0}

In this section, we will prove the existence and uniqueness of problem
\eqref{1.1}-\eqref{1.3} by using Faedo-Galerkin method.

\noindent{\bf Proof of Theorem \ref{th2.1}.}
We divide the proof of Theorem \ref{th2.1} in two steps.

Step 1: Faedo-Galerkin approximation.

We construct approximations of the solution $(u,v,z)$ by the Faedo-Galerkin method as follows. For $n\geq 1$, let
$W_{n}={\rm span}\{w_{1}, \ldots, w_{i}\}$ be a Hilbertian basis of the space $H^{1}(\Omega)$ and $Y_{n}={\rm span}\{\psi_{1}, \ldots, \psi_{i}\}$ be a Hilbertian basis of the space $H^{1}(L_{1},L_{2})$.

Now, we define for $1\leq j\leq n$ the sequence $\varphi_{j}(x,\rho)$ as follows:
\begin{equation*}
\varphi_{j}(x,0)=w_{j}(x).
\end{equation*}
Then we may extend $\varphi_{j}(x,0)$ by $\varphi_{j}(x,\rho)$ over $L^2((0,1),H^{1}(\Omega))$ and denote
$V_{n}={\rm span}\{\varphi_{1}, \ldots, \varphi_{n}\}$.

We choose sequences $\left(u_{0}^{(n)}\right)$, $\left(u_{1}^{(n)}\right)$ in $W_{n}$, $\left(v_{0}^{(n)}\right)$, $\left(v_{1}^{(n)}\right)$ in $Y_{n}$ and a sequence $\left(z_{0}^{(n)}\right)$ in $V_{n}$
such that $u_{0}^{(n)}\rightarrow u_{0}$ strongly in $H^{1}(\Omega)$, $u_{1}^{(n)}\rightarrow u_{1}$ strongly in $H^{1}(\Omega)$,
$v_{0}^{(n)}\rightarrow v_{0}$ strongly in $H^{1}(L_{1},L_{2})$, $v_{1}^{(n)}\rightarrow v_{1}$ strongly in $H^{1}(L_{1},L_{2})$ and
$z_{0}^{(n)}\rightarrow f_{0}$ strongly in $L^2((0,1),H^{1}(\Omega))$.

We define the approximations
\begin{equation*}
\left(u^{(n)}(x,t),v^{(n)}(x,t)\right)=\sum\limits_{i=1}^{n}h_{i}^{(n)}(t)(w_{i}(x),\psi_{i}(x))\quad {\rm and} \quad z^{(n)}(x,\rho,t)=\sum\limits_{i=1}^{n}f_{i}^{(n)}(t)\varphi_{i}(x),
\end{equation*}
where $\left(u^{(n)}(t),v^{(n)}(t),z^{(n)}(t)\right)$ is a solution to the following Cauchy problem:
\bes\left\{\begin{array}{ll}
\displaystyle \int_{\Omega}u_{tt}^{(n)}w_{i}{\rm d}x-\left[\left(au_{x}^{(n)}-g*u_{x}^{(n)}\right)w_{i}\right]_{\partial\Omega}+\int_{\Omega}au_{x}^{(n)}w_{ix}{\rm d}x-\int_{\Omega}\left(g*u_{x}^{(n)}\right)w_{ix}{\rm d}x\\
 \displaystyle+\int_{\Omega}\mu_{1}u_{t}^{(n)}w_{i}{\rm d}x+\int_{\Omega}\mu_{2}z^{(n)}(x,1,t)w_{i}{\rm d}x=0,\medskip
\\
  \displaystyle \int_{L_{1}}^{L_{2}}v_{tt}^{(n)}\psi_{i}{\rm d}x+\int_{L_{1}}^{L_{2}}bv_{x}^{(n)}\psi_{ix}{\rm d}x-\left[bv_{x}^{(n)}\psi_{i}\right]_{L_{1}}^{L_{2}}=0, \medskip\\
\displaystyle  z^{(n)}(x,0,t)=u_{t}^{(n)}(x,t),
  \medskip\\
\displaystyle \left(u^{(n)}(0),u_{t}^{(n)}(0)\right)=\left(u_{0}^{(n)},u_{1}^{(n)}\right)
 \end{array}\right.\label{2.21}
 \ees
and
\bes\left\{\begin{array}{ll}
\displaystyle \int_{\Omega}\left(\tau z_{t}^{(n)}(x,\rho,t)+z_{\rho}^{(n)}(x,\rho,t)\right)\varphi_{i}{\rm d}x=0,
\\
  \displaystyle z^{(n)}(\rho,0)=z_{0}^{(n)}.
 \end{array}\right.\label{2.22}
 \ees
According to the standard theory of ordinary differential equations, the finite dimensional problem \eqref{2.21}-\eqref{2.22} have a
solution $\left(h_{i}^{(n)}(t),f_{i}^{(n)}(t)\right)_{i=1, \ldots, n}$ defined on $[0,t_{n})$.

Step 2: Energy estimates.

Multiplying the first and the second equation of \eqref{2.21} by $\left(h_{i}^{(n)}\right)'(t)$, we have
\begin{align}\label{2.23}
&\int_{\Omega}u_{tt}^{(n)}u_{t}^{(n)}{\rm d}x-\left[\left(au_{x}^{(n)}-g*u_{x}^{(n)}\right)w_{i}\right]_{\partial\Omega}\times \left(h_{i}^{(n)}\right)'(t)+\int_{\Omega}au_{x}^{(n)}u_{xt}^{(n)}{\rm d}x\nonumber\\
&-\int_{\Omega}\left(g*u_{x}^{(n)}\right)u_{xt}^{(n)}{\rm d}x+\int_{\Omega}\mu_{1}u_{t}^{(n)}u_{t}^{(n)}{\rm d}x+\int_{\Omega}\mu_{2}z^{(n)}(x,1,t)u_{t}^{(n)}{\rm d}x=0
\end{align}
and
\begin{equation}\label{2.24}
\int_{L_{1}}^{L_{2}}v_{tt}^{(n)}v_{t}^{(n)}{\rm d}x+\int_{L_{1}}^{L_{2}}bv_{x}^{(n)}v_{xt}^{(n)}{\rm d}x-\left[bv_{x}^{(n)}\psi_{i}\right]_{L_{1}}^{L_{2}}\times \left(h_{i}^{(n)}\right)'(t)=0.
\end{equation}
Multiplying the first equation of \eqref{2.22} by $\frac{\zeta}{\tau}f_{i}^{(n)}(t)$ and
integrating over $(0,t)\times(0,1)$, we get
\begin{align}\label{2.25}
&\frac{\zeta}{2}\int_{\Omega}\int_{0}^{1}\left(z^{(n)}\right)^2(x,\rho,t){\rm d}\rho{\rm d}x+\frac{\zeta}{\tau}\int_{0}^{t}\int_{\Omega}\int_{0}^{1}z_{\rho}^{(n)}z^{(n)}(x,\rho,s){\rm d}\rho{\rm d}x{\rm d}s\nonumber\\
=&\frac{\zeta}{2}\int_{\Omega}\int_{0}^{1}\left(z_{0}^{(n)}\right)^2{\rm d}\rho{\rm d}x.
\end{align}
To handle the last term in the left-hand side of \eqref{2.25}, we remark that
\begin{align}\label{2.26}
\int_{0}^{t}\int_{\Omega}\int_{0}^{1}z_{\rho}^{(n)}z^{(n)}(x,\rho,s){\rm d}\rho{\rm d}x{\rm d}s&=\frac{1}{2}\int_{0}^{t}\int_{\Omega}\int_{0}^{1}\frac{\partial}{\partial\rho}\left(z^{(n)}\right)^2(x,\rho,s){\rm d}\rho{\rm d}x{\rm d}s\nonumber\\
&=\frac{1}{2}\int_{0}^{t}\int_{\Omega}\left(\left(z^{(n)}\right)^2(x,1,s)-\left(z^{(n)}\right)^2(x,0,s)\right){\rm d}x{\rm d}s.
\end{align}
Integrating \eqref{2.23} and \eqref{2.24} over $(0,t)$, counting them and \eqref{2.25} up, taking into account \eqref{2.26}
and using Lemma \ref{le2.1}, we obtain
\begin{align}\label{2.27}
&\mathscr{E}_{n}(t)+\left(\mu_{1}-\frac{\zeta}{2\tau}\right)\int_{0}^{t}\int_{\Omega}\left(u_{t}^{(n)}\right)^2(x,s){\rm d}x{\rm d}s
+\frac{\zeta}{2\tau}\int_{0}^{t}\int_{\Omega}\left(z^{(n)}\right)^2(x,1,s){\rm d}x{\rm d}s\nonumber\\
&+\mu_{2}\int_{0}^{t}\int_{\Omega}z^{(n)}(x,1,s)u_{t}^{(n)}(x,s){\rm d}x{\rm d}s+\frac{1}{2}\int_{0}^{t}\int_{\Omega}g(t)\left|u_{x}^{(n)}\right|^2{\rm d}x{\rm d}s-\frac{1}{2}\int_{0}^{t}\int_{\Omega}\left(g'\square u_{x}^{(n)}\right){\rm d}x{\rm d}s\nonumber\\
=&\mathscr{E}_{n}(0),
\end{align}
where
\begin{align}\label{2.28}
\mathscr{E}_{n}(t)=&\frac{1}{2}\int_{\Omega}\left(u_{t}^{(n)}\right)^2(x,t){\rm d}x+\frac{1}{2}\beta(t)\int_{\Omega}\left(u_{x}^{(n)}\right)^2(x,t){\rm d}x+\frac{1}{2}\int_{\Omega}\left(g\square u_{x}^{(n)}\right){\rm d}x\nonumber\\
&+\frac{1}{2}\int_{L_{1}}^{L_{2}}\left[\left(v_{t}^{(n)}\right)^2(x,t)+b\left(v_{x}^{(n)}\right)^2(x,t)\right]{\rm d}x+\frac{\zeta}{2}\int_{\Omega}\int_{0}^{1}\left(z^{(n)}\right)^2(x,\rho,t){\rm d}\rho{\rm d}x.
\end{align}
At this point, we have to distinguish the following two cases:

Case 1: We suppose that $\mu_{2}<\mu_{1}$ and choose $\zeta$ satisfying \eqref{2.9}. Young's inequality gives us that
\begin{align*}
&\mathscr{E}_{n}(t)+\left(\mu_{1}-\frac{\zeta}{2\tau}-\frac{\mu_{2}}{2}\right)\int_{0}^{t}\int_{\Omega}\left(u_{t}^{(n)}\right)^2(x,s){\rm d}x{\rm d}s
+\left(\frac{\zeta}{2\tau}-\frac{\mu_{2}}{2}\right)\int_{0}^{t}\int_{\Omega}\left(z^{(n)}\right)^2(x,1,s){\rm d}x{\rm d}s\nonumber\\
&+\frac{1}{2}\int_{0}^{t}\int_{\Omega}g(t)\left|u_{x}^{(n)}\right|^2{\rm d}x{\rm d}s-\frac{1}{2}\int_{0}^{t}\int_{\Omega}\left(g'\square u_{x}^{(n)}\right){\rm d}x{\rm d}s\nonumber\\
\leq&\mathscr{E}_{n}(0).
\end{align*}
Consequently, using \eqref{2.9}, we have
\begin{align}\label{2.29}
&\mathscr{E}_{n}(t)+c_{1}\int_{0}^{t}\int_{\Omega}\left(u_{t}^{(n)}\right)^2(x,s){\rm d}x{\rm d}s
+c_{2}\int_{0}^{t}\int_{\Omega}\left(z^{(n)}\right)^2(x,1,s){\rm d}x{\rm d}s
\nonumber\\
&+\frac{1}{2}\int_{0}^{t}\int_{\Omega}g(t)\left|u_{x}^{(n)}\right|^2{\rm d}x{\rm d}s-\frac{1}{2}\int_{0}^{t}\int_{\Omega}\left(g'\square u_{x}^{(n)}\right){\rm d}x{\rm d}s\nonumber\\
\leq&\mathscr{E}_{n}(0).
\end{align}

Case 2: We suppose that $\mu_{2}=\mu_{1}=\mu$ and choose $\zeta=\tau\mu$. Then \eqref{2.29} takes the form
\begin{equation}
\mathscr{E}_{n}(t)+\frac{1}{2}\int_{0}^{t}\int_{\Omega}g(t)\left|u_{x}^{(n)}\right|^2{\rm d}x{\rm d}s-\frac{1}{2}\int_{0}^{t}\int_{\Omega}\left(g'\square u_{x}^{(n)}\right){\rm d}x{\rm d}s\leq\mathscr{E}_{n}(0).
\end{equation}

Now, since the sequences $\left(u_{0}^{(n)}\right)_{n\in \mathbb{N}}$, $\left(u_{1}^{(n)}\right)_{n\in \mathbb{N}}$, $\left(v_{0}^{(n)}\right)_{n\in \mathbb{N}}$, $\left(v_{1}^{(n)}\right)_{n\in \mathbb{N}}$, $\left(z_{0}^{(n)}\right)_{n\in \mathbb{N}}$ converge and using (G2),
in the both cases we can find a positive constant $c_{3}$ independent of $n$ such that
\begin{equation}\label{2.30}
\mathscr{E}_{n}(t)\leq c_{3}.
\end{equation}
Therefore, from \eqref{2.30} and the Lion-Aubin's compactness theorem \cite{l1969}, we can pass to the limit in \eqref{2.21}.
The rest of the proof is routine.

\section{General decay of the solution}\label{s4}
\setcounter{equation}{0}

In this section, we consider the asymptotic behavior of problem \eqref{1.1}-\eqref{1.3}.
For the proof of Theorem \ref{th3.1}, we use the following lemmas.
\begin{lemma}
Let $(u,v,z)$ be the solution of problem \eqref{2.5}-\eqref{2.6}. Assume that $\mu_{2}<\mu_{1}$. Then we have the inequality
\begin{equation}\label{3.61}
\frac{d}{dt}E(t)\leq-c_{4}\left[\int_{\Omega}u_{t}^2(x,t){\rm d}x+\int_{\Omega}z^2(x,1,t){\rm d}x\right]+\frac{1}{2}\int_{\Omega}(g'\square u_{x})(t){\rm d}x.
\end{equation}
\end{lemma}
{\bf Proof.}
Multiplying the first equation of \eqref{2.5} by $u_{t}$, the second equation of \eqref{2.5} by $v_{t}$,
integrating by parts and \eqref{2.6}, we obtain
\begin{align}\label{3.2}
&\frac{1}{2}\frac{d}{dt}\left\{\int_{\Omega}[u_{t}^2(x,t)+au_{x}^2(x,t)]{\rm d}x\right\}+\frac{1}{2}\frac{d}{dt}\left\{\int_{L_{1}}^{L_{2}}[v_{t}^2(x,t)+bv_{x}^2(x,t)]{\rm d}x\right\}\nonumber\\
=&-\mu_{1}\int_{\Omega}u_{t}^2(x,t){\rm d}x-\mu_{2}\int_{\Omega}u_{t}(x,t)z(x,1,t){\rm d}x+\int_{0}^{t}g(t-s)\int_{\Omega}u_{x}(s)u_{xt}(t){\rm d}s{\rm d}x.
\end{align}
From Lemma \ref{le2.1}, the last term in the right-hand side of \eqref{3.2} can be rewritten as
\begin{align*}
&\int_{0}^{t}g(t-s)\int_{\Omega}u_{x}(s)u_{xt}(t){\rm d}s{\rm d}x+\frac{1}{2}g(t)\int_{\Omega}u_{x}^2{\rm d}x\\
=&\frac{1}{2}\frac{d}{dt}\left\{\int_{0}^{t}g(s)\int_{\Omega}u_{x}^2{\rm d}x{\rm d}s-\int_{\Omega}(g\square u_{x})(t){\rm d}x\right\}+\frac{1}{2}\int_{\Omega}(g'\square u_{x})(t){\rm d}x.
\end{align*}
So \eqref{3.2} becomes
\begin{align}\label{3.2'}
&\frac{1}{2}\frac{d}{dt}\left\{\int_{\Omega}\left[u_{t}^2(x,t)+\beta(t)u_{x}^2(x,t)\right]{\rm d}x\right\}+\frac{1}{2}\frac{d}{dt}\left\{\int_{L_{1}}^{L_{2}}[v_{t}^2(x,t)+bv_{x}^2(x,t)]{\rm d}x\right\}\nonumber\\
&+\frac{1}{2}\frac{d}{dt}\int_{\Omega}(g\square u_{x})(t){\rm d}x\nonumber\\
=&-\mu_{1}\int_{\Omega}u_{t}^2(x,t){\rm d}x-\mu_{2}\int_{\Omega}u_{t}(x,t)z(x,1,t){\rm d}x-\frac{1}{2}g(t)\int_{\Omega}u_{x}^2{\rm d}x+\frac{1}{2}\int_{\Omega}(g'\square u_{x})(t){\rm d}x.
\end{align}

Now, multiplying the third equation of \eqref{2.5} by $\frac{\zeta}{\tau}z$ and integrating the result over $\Omega\times(0,1)$ with
respect to $x$ and $\rho$ respectively, we have
\begin{equation}\label{3.3}
\frac{\zeta}{2}\frac{d}{dt}\int_{\Omega}\int_{0}^{1}z^2(x,\rho,t){\rm d}\rho{\rm d}x=-\frac{\zeta}{2\tau}\int_{\Omega}(z^2(x,1)-z^2(x,0)){\rm d}x.
\end{equation}
Using \eqref{3.1}, \eqref{3.2'} and \eqref{3.3}, we gain
\begin{align}\label{3.4}
\frac{d}{dt}E(t)=&-\left(\mu_{1}-\frac{\zeta}{2\tau}\right)\int_{\Omega}u_{t}^2(x,t){\rm d}x-\frac{\zeta}{2\tau}\int_{\Omega}z^2(x,1,t){\rm d}x-\mu_{2}\int_{\Omega}u_{t}(x,t)z(x,1,t){\rm d}x\nonumber\\
&-\frac{1}{2}g(t)\int_{\Omega}u_{x}^2{\rm d}x+\frac{1}{2}\int_{\Omega}(g'\square u_{x})(t){\rm d}x.
\end{align}
By Young's inequality in \eqref{3.4}, we get
\begin{align*}
\frac{d}{dt}E(t)\leq&-\left(\mu_{1}-\frac{\zeta}{2\tau}-\frac{\mu_{2}}{2}\right)\int_{\Omega}u_{t}^2(x,t){\rm d}x-\left(\frac{\zeta}{2\tau}-\frac{\mu_{2}}{2}\right)\int_{\Omega}z^2(x,1,t){\rm d}x\nonumber\\
&+\frac{1}{2}\int_{\Omega}(g'\square u_{x})(t){\rm d}x.
\end{align*}
Then exploiting \eqref{2.9} our conclusion holds. The proof is complete.

Now, we define the functional $\mathscr{D}(t)$ as follows
\begin{equation*}
\mathscr{D}(t)=\int_{\Omega}uu_{t}{\rm d}x+\frac{\mu_{1}}{2}\int_{\Omega}u^2{\rm d}x+\int_{L_{1}}^{L_{2}}vv_{t}{\rm d}x.
\end{equation*}
Then we have the following estimate.
\begin{lemma}
The functional $\mathscr{D}(t)$ satisfies
\begin{align}\label{3.5}
\frac{d}{dt}\mathscr{D}(t)\leq&\int_{\Omega}u_{t}^2{\rm d}x+\int_{L_{1}}^{L_{2}}v_{t}^2{\rm d}x+(c^{*}\varepsilon +\varepsilon-\beta(t))\int_{\Omega}u_{x}^2{\rm d}x+\frac{1}{4\varepsilon}(a-\beta(t))\int_{\Omega}(g\square u_{x}){\rm d}x\nonumber\\
&+\frac{\mu_{2}^2}{4\varepsilon}\int_{\Omega}z^2(x,1,t){\rm d}x-\int_{L_{1}}^{L_{2}}bv_{x}^2{\rm d}x.
\end{align}
\end{lemma}
{\bf Proof.}
Taking the derivative of $\mathscr{D}(t)$ with respect to $t$ and using \eqref{2.5}, we have
\begin{align}\label{3.6}
\frac{d}{dt}\mathscr{D}(t)=&\int_{\Omega}u_{t}^2{\rm d}x-\int_{\Omega}(au_{x}-g*u_{x})u_{x}{\rm d}x-\mu_{2}\int_{\Omega}z(x,1,t)u{\rm d}x
+\int_{L_{1}}^{L_{2}}v_{t}^2{\rm d}x-\int_{L_{1}}^{L_{2}}bv_{x}^2{\rm d}x\nonumber\\
=&\int_{\Omega}u_{t}^2{\rm d}x-\beta(t)\int_{\Omega}u_{x}^2{\rm d}x-\int_{\Omega}(g\diamond u_{x})u_{x}{\rm d}x-\mu_{2}\int_{\Omega}z(x,1,t)u{\rm d}x
+\int_{L_{1}}^{L_{2}}v_{t}^2{\rm d}x\nonumber\\
&-\int_{L_{1}}^{L_{2}}bv_{x}^2{\rm d}x.
\end{align}
By exploiting Young's inequality and Poincar$\acute{e}$'s inequality, we get for any $\varepsilon>0$
\begin{equation}\label{3.7}
\mu_{2}\int_{\Omega}z(x,1,t)u{\rm d}x\leq\frac{\mu_{2}^2}{4\varepsilon}\int_{\Omega}z^2(x,1,t){\rm d}x+c^{*}\varepsilon\int_{\Omega}u_{x}^2{\rm d}x.
\end{equation}
Young's inequality and (G1) imply that
\begin{align}\label{3.8}
\int_{\Omega}(g\diamond u_{x})u_{x}{\rm d}x&\leq\varepsilon\int_{\Omega}u_{x}^2{\rm d}x+\frac{1}{4\varepsilon}\int_{\Omega}(g\diamond u_{x})^2{\rm d}x\nonumber\\
&\leq\varepsilon\int_{\Omega}u_{x}^2{\rm d}x+\frac{1}{4\varepsilon}(a-\beta(t))\int_{\Omega}(g\square u_{x}){\rm d}x.
\end{align}
Inserting the estimates \eqref{3.7} and \eqref{3.8} into \eqref{3.6}, then \eqref{3.5} is fulfilled. The proof is complete.

Now, inspired by \cite{mmn2002}, we introduce the function
\begin{align}
q(x)=\left\{{\begin{aligned} &x-\frac{L_{1}}{2},&&x\in[0,L_{1}],\\
&\frac{L_{1}}{2}-\frac{L_{1}+L_{3}-L_{2}}{2(L_{2}-L_{1})}(x-L_{1}),&&x\in(L_{1},L_{2}),\\
&x-\frac{L_{2}+L_{3}}{2},&&x\in[L_{2},L_{3}].\\
\end{aligned}}\right.
\end{align}
It is easy to see that $q(x)$ is bounded, that is $|q(x)|\leq M$, where $M=\max\left\{\frac{L_{1}}{2},\frac{L_{3}-L_{2}}{2}\right\}$ is a positive constant. And we define the functionals
\begin{equation}
\mathscr{F}_{1}(t)=-\int_{\Omega}q(x)u_{t}(au_{x}-g*u_{x}){\rm d}x,\quad \mathscr{F}_{2}(t)=-\int_{L_{1}}^{L_{2}}q(x)v_{x}v_{t}{\rm d}x.
\end{equation}
Then we have the following estimates.
\begin{lemma}\label{le3.5}
The functionals $\mathscr{F}_{1}(t)$ and $\mathscr{F}_{2}(t)$ satisfy
\begin{align}\label{3.9}
\frac{d}{dt}\mathscr{F}_{1}(t)\leq&\left[-\frac{q(x)}{2}(au_{x}-g*u_{x})^2\right]_{\partial\Omega}-\left[\frac{a}{2}q(x)u_{t}^2\right]_{\partial\Omega}
+\left[\frac{a}{2}+\frac{\mu_{1}^2}{2\varepsilon_{1}}+\frac{M^2}{4\varepsilon_{1}}\right]\int_{\Omega}u_{t}^2{\rm d}x\nonumber\\
&+\left[\varepsilon_{1}M^2a^2+\beta^2(t)+2M^2\varepsilon_{1}(a-\beta(t))^2+c_{5}^2\varepsilon_{1}\right]\int_{\Omega}u_{x}^2{\rm d}x+\frac{\mu_{2}^2}{2\varepsilon_{1}}\int_{\Omega}z^2(x,1,t){\rm d}x\nonumber\\
&+(1+2 M^2\varepsilon_{1})(a-\beta(t))\int_{\Omega}(g\square u_{x}){\rm d}x+(a-\beta(t))\varepsilon_{1}\int_{\Omega}(g'\square u_{x}){\rm d}x
\end{align}
and
\begin{align}\label{3.10}
\frac{d}{dt}\mathscr{F}_{2}(t)\leq&-\frac{L_{1}+L_{3}-L_{2}}{4(L_{2}-L_{1})}\left(\int_{L_{1}}^{L_{2}}v_{t}^2{\rm d}x+\int_{L_{1}}^{L_{2}}bv_{x}^2{\rm d}x\right)+\frac{L_{1}}{4}v_{t}^2(L_{1})+\frac{L_{3}-L_{2}}{4}v_{t}^2(L_{2})\nonumber\\
&+\frac{b}{4}\left((L_{3}-L_{2})v_{x}^2(L_{2},t)+L_{1}v_{x}^2(L_{1},t)\right).
\end{align}
\end{lemma}
{\bf Proof.}
Taking the derivative of $\mathscr{F}_{1}(t)$ with respect to $t$ and using \eqref{2.5}, we get
\begin{align}\label{3.11}
\frac{d}{dt}\mathscr{F}_{1}(t)=&-\int_{\Omega}q(x)u_{tt}(au_{x}-g*u_{x}){\rm d}x-\int_{\Omega}q(x)u_{t}\left(au_{xt}-g(t)u_{x}(t)+(g'\diamond u_{x})(t)\right){\rm d}x\nonumber\\
=&\left[-\frac{q(x)}{2}(au_{x}-g*u_{x})^2\right]_{\partial\Omega}+\frac{1}{2}\int_{\Omega}q'(x)(au_{x}-g*u_{x})^2{\rm d}x-\left[\frac{a}{2}q(x)u_{t}^2\right]_{\partial\Omega}\nonumber\\
&+\frac{a}{2}\int_{\Omega}q'(x)u_{t}^2{\rm d}x-\int_{\Omega}q(x)(\mu_{1}u_{t}(x,t)+\mu_{2}z(x,1,t))(g*u_{x}){\rm d}x\nonumber\\
&+\int_{\Omega}q(x)au_{x}(\mu_{1}u_{t}(x,t)+\mu_{2}z(x,1,t)){\rm d}x-\int_{\Omega}q(x)u_{t}[(g'\diamond u_{x})(t)-g(t)u_{x}]{\rm d}x.
\end{align}
We  note that
\begin{align}\label{3.12}
\frac{1}{2}\int_{\Omega}q'(x)(au_{x}-g*u_{x})^2{\rm d}x&=\frac{1}{2}\int_{\Omega}\left[\left(a-\int_{0}^{t}g(s){\rm d}s\right)u_{x}+g\diamond u_{x}\right]^2{\rm d}x\nonumber\\
&\leq\int_{\Omega}|\beta(t)|^2u_{x}^2{\rm d}x+\int_{\Omega}|g\diamond u_{x}|^2{\rm d}x\nonumber\\
&\leq\int_{\Omega}|\beta(t)|^2u_{x}^2{\rm d}x+(a-\beta(t))\int_{\Omega}(g\square u_{x}){\rm d}x.
\end{align}
Young's inequality gives us  for any $\varepsilon_{1}>0$,
\begin{align}\label{3.13}
&\int_{\Omega}q(x)au_{x}(\mu_{1}u_{t}(x,t)+\mu_{2}z(x,1,t)){\rm d}x\nonumber\\
\leq&\varepsilon_{2}M^2 a^2\int_{\Omega}u_{x}^2{\rm d}x+\frac{\mu_{1}^2}{4\varepsilon_{1}}\int_{\Omega}u_{t}^2{\rm d}x+\frac{\mu_{2}^2}{4\varepsilon_{1}}\int_{\Omega}z^2(x,1,t){\rm d}x,
\end{align}
\begin{align}\label{3.14}
&\int_{\Omega}q(x)(\mu_{1}u_{t}(x,t)+\mu_{2}z(x,1,t))(g*u_{x}){\rm d}x\nonumber\\
\leq&\varepsilon_{1}M^2\int_{\Omega}(g*u_{x})^2{\rm d}x+\frac{\mu_{1}^2}{4\varepsilon_{1}}\int_{\Omega}u_{t}^2{\rm d}x+\frac{\mu_{2}^2}{4\varepsilon_{1}}\int_{\Omega}z^2(x,1,t){\rm d}x\nonumber\\
\leq&2\varepsilon_{1}M^2(a-\beta(t))^2\int_{\Omega}u_{x}^2{\rm d}x+2M^2\varepsilon_{1}(a-\beta(t))\int_{\Omega}(g\square u_{x}){\rm d}x+\frac{\mu_{1}^2}{4\varepsilon_{1}}\int_{\Omega}u_{t}^2{\rm d}x\nonumber\\
&+\frac{\mu_{2}^2}{4\varepsilon_{1}}\int_{\Omega}z^2(x,1,t){\rm d}x
\end{align}
and
\begin{align}\label{3.15}
&\int_{\Omega}q(x)u_{t}[(g'\diamond u_{x})(t)-g(t)u_{x}]{\rm d}x\nonumber\\
\leq&\frac{M^2}{4\varepsilon_{1}}\int_{\Omega}u_{t}^2{\rm d}x+c_{5}^2\varepsilon_{1}\int_{\Omega}u_{x}^2{\rm d}x+(a-\beta(t))\varepsilon_{1}\int_{\Omega}(g'\square u_{x}){\rm d}x.
\end{align}
Inserting \eqref{3.12}-\eqref{3.15} into \eqref{3.11}, we get \eqref{3.9}.

By the same method, taking the derivative of $\mathscr{F}_{1}(t)$ with respect to $t$, we obtain
\begin{align*}
\frac{d}{dt}\mathscr{F}_{2}(t)=&-\int_{L_{1}}^{L_{2}}q(x)v_{xt}v_{t}{\rm d}x-\int_{L_{1}}^{L_{2}}q(x)v_{x}v_{tt}{\rm d}x\\
=&\left[-\frac{1}{2}q(x)v_{t}^2\right]_{L_{1}}^{L_{2}}+\frac{1}{2}\int_{L_{1}}^{L_{2}}q'(x)v_{t}^2{\rm d}x+\frac{1}{2}\int_{L_{1}}^{L_{2}}bq'(x)v_{x}^2{\rm d}x+\left[-\frac{b}{2}q(x)v_{x}^2\right]_{L_{1}}^{L_{2}}\\
\leq&-\frac{L_{1}+L_{3}-L_{2}}{4(L_{2}-L_{1})}\left(\int_{L_{1}}^{L_{2}}v_{t}^2{\rm d}x+\int_{L_{1}}^{L_{2}}bv_{x}^2{\rm d}x\right)+\frac{L_{1}}{4}v_{t}^2(L_{1})+\frac{L_{3}-L_{2}}{4}v_{t}^2(L_{2})\\
&+\frac{b}{4}\left((L_{3}-L_{2})v_{x}^2(L_{2},t)+L_{1}v_{x}^2(L_{1},t)\right).
\end{align*}
Thus, the proof of Lemma \ref{le3.5} is finished.

As in \cite{amm2012}, we define the functional
\begin{equation*}
\mathscr{F}_{3}(t)=\tau\int_{\Omega}\int_{0}^{1}e^{-\tau\rho}z^2(x,\rho,t){\rm d}\rho{\rm d}x,
\end{equation*}
then we have the following estimate.
\begin{lemma}{\rm (\cite{amm2012})}
The functionals $\mathscr{F}_{3}(t)$ satisfies
\begin{equation*}
\frac{d}{dt}\mathscr{F}_{3}(t)\leq-c_{6}\left(\int_{\Omega}z^2(x,1,t){\rm d}x+\tau\int_{\Omega}\int_{0}^{1}z^2(x,\rho,t){\rm d}\rho{\rm d}x\right)+\int_{\Omega}u_{t}^2(x,t){\rm d}x.
\end{equation*}
\end{lemma}

Now, we are ready to prove Theorem \ref{th3.1}.

\noindent{\bf Proof of Theorem \ref{th3.1}. } We define the Lyapunov functional
\begin{equation}\label{3.20}
L(t)=N_{1}E(t)+N_{2}\mathscr{D}(t)+N_{3}\mathscr{F}_{1}(t)+N_{4}\mathscr{F}_{2}(t)+\mathscr{F}_{3}(t),
\end{equation}
where $N_{1}, N_{2}, N_{3}$ and $N_{4}$ are positive constants that will be fixed later.

Taking the derivative of \eqref{3.20} with respect to $t$ and making the use of the above lemmas, we have
\begin{align}\label{3.21}
\frac{d}{dt}L(t)\leq&\left\{-N_{1}c_{4}+1+N_{2}+N_{3}\left(\frac{a}{2}+\frac{\mu_{1}^2}{2\varepsilon_{1}}+\frac{M^2}{4\varepsilon_{1}}\right)\right\}
\int_{\Omega}u_{t}^2{\rm d}x\nonumber\\
&+\left\{-N_{1}c_{4}-c_{6}+\frac{\mu_{2}^2N_{2}}{4\varepsilon}+\frac{\mu_{2}^2 N_{3}}{2\varepsilon_{1}}\right\}
\int_{\Omega}z^2(x,1,t){\rm d}x\nonumber\\
&+\left\{-N_{2}(\beta(t)-c^{*}\varepsilon-\varepsilon)+N_{3}\left(\varepsilon_{1}M^2a^2+\beta(t)^2+2M^2\varepsilon_{1}(a-\beta(t))^2
+c_{5}^2\varepsilon_{1}\right)\right\}\int_{\Omega}u_{x}^2{\rm d}x\nonumber\\
&+\left\{-\frac{b(L_{1}+L_{3}-L_{2})}{4(L_{2}-L_{1})}N_{4}-N_{2}b\right\}\int_{L_{1}}^{L_{2}}v_{x}^2{\rm d}x\nonumber\\
&+\left\{-\frac{L_{1}+L_{3}-L_{2}}{4(L_{2}-L_{1})}N_{4}+N_{2}\right\}\int_{L_{1}}^{L_{2}}v_{t}^2{\rm d}x\nonumber\\
&+(N_{4}-bN_{3})\frac{b}{4}\left((L_{3}-L_{2})v_{x}^2(L_{2},t)+L_{1}v_{x}^2(L_{1},t)\right)\nonumber\\
&+(N_{4}-aN_{3})\left[\frac{L_{1}}{4}v_{t}^2(L_{1},t)
+\frac{L_{3}-L_{2}}{4}v_{t}^2(L_{2},t)\right]\nonumber\\
&+c(N_{2},N_{3})\int_{\Omega}(g\square u_{x}){\rm d}x+\left(\frac{N_{1}}{2}-c(N_{3})\right)
\int_{\Omega}(g'\square u_{x}){\rm d}x.
\end{align}
At this moment, we wish all coefficients except the last two in \eqref{3.21} will be negative. In fact, under assumption \eqref{3.30}, we can find $N_{2}, N_{3}$ and $N_{4}$ such that
\begin{align*}
N_{2}<\frac{L_{1}+L_{3}-L_{2}}{4(L_{2}-L_{1})}N_{4},\quad N_{4}<bN_{3},\quad N_{4}<aN_{3},\quad N_{2}>N_{3}\beta_{0}.
\end{align*}
Once the above constants are fixed, we may choose $\varepsilon$ and $\varepsilon_{1}$ small enough such that
\begin{equation*}
N_{2}(c^{*}\varepsilon+\varepsilon)+N_{3}\left(\varepsilon_{1}M^2a^2+2M^2\varepsilon_{1}(a-\beta(t))^2
+c_{5}^2\varepsilon_{1}\right)<N_{2}-N_{3}\beta(t).
\end{equation*}
Finally, choosing $N_{1}$ large enough such that the first two coefficients in \eqref{3.21} are negative and the last coefficient is positive.
From the above, we deduce that, there exists two positive constants $\alpha_{1}$ and $\alpha_{2}$ such that \eqref{3.21} becomes
\begin{equation}\label{3.41}
\frac{d}{dt}L(t)\leq-\alpha_{1}E(t)+\alpha_{2}\int_{\Omega}(g\square u_{x}){\rm d}x.
\end{equation}
On the other hand, by the definition of the functionals $\mathscr{D}(t)$, $\mathscr{F}_{1}(t)$, $\mathscr{F}_{2}(t)$, $\mathscr{F}_{3}(t)$ and $E(t)$, for $N_{1}$ large enough, there exists a positive constant $\alpha_{3}$ satisfying
\begin{equation*}
|N_{2}\mathscr{D}(t)+N_{3}\mathscr{F}_{1}(t)+N_{4}\mathscr{F}_{2}(t)+\mathscr{F}_{3}(t)|\leq\alpha_{3} E(t),
\end{equation*}
which implies that
\begin{equation*}
(N_{1}-\alpha_{3})E(t)\leq L(t)\leq(N_{1}+\alpha_{3})E(t).
\end{equation*}

In order to finish the proof of the stability estimates, we need to estimate the last term in \eqref{3.41}.
Exploiting (G2) and \eqref{3.61}, we have
\begin{align}\label{3.42}
 \xi(t)\int_{\Omega}(g\square u_{x}){\rm d}x
\leq \int_{\Omega}[(\xi g)\square u_{x}]{\rm d}x
\leq -\int_{\Omega}(g'\square u_{x}){\rm d}x
\leq -2\frac{d}{dt}E(t).
\end{align}

Now, we define functionals $\mathscr{L}(t)$ as
\begin{equation*}
\mathscr{L}(t)=\xi(t)L(t)+2\alpha_{2}E(t).
\end{equation*}
The fact that $L(t)$ and $E(t)$ are equivalent and (G2) imply that, for some positive constants $\eta_{1}$ and $\eta_{2}$,
\begin{equation}\label{3.44}
\eta_{1}E(t)\leq\mathscr{L}(t)\leq\eta_{2}E(t),
\end{equation}
Using \eqref{3.42}, \eqref{3.44} and (G2), we obtain
\begin{align*}
\frac{d}{dt}\mathscr{L}(t)&=\xi'(t)L(t)+\xi(t)\frac{d}{dt}L(t)+2\alpha_{2}\frac{d}{dt}E(t)\\
&\leq\xi(t)\left(-\alpha_{1}E(t)+\alpha_{2}\int_{\Omega}(g\square u_{x}){\rm d}x\right)+2\alpha_{2}\frac{d}{dt}E(t)\\
&\leq-\alpha_{1}\xi(t)E(t)\\
&\leq-\gamma_{0}\xi(t)\mathscr{L}(t),
\end{align*}
where $\gamma_{0}=\frac{\alpha_{1}}{\eta_{2}}$.
We conclude that, for any $\gamma_{1}\in(0,\gamma_{0})$,
\begin{equation}\label{3.45}
\frac{d}{dt}\mathscr{L}(t)\leq-\gamma_{1}\xi(t)\mathscr{L}(t).
\end{equation}
A simple integration of \eqref{3.45} leads to
\begin{equation}\label{3.25}
\mathscr{L}(t)\leq\mathscr{L}(0)e^{-\gamma_{1}\int_{0}^{t}\xi(s){\rm d}s},\quad \forall t\geq0.
\end{equation}
Again, use of \eqref{3.44} and \eqref{3.25} yields the desired result \eqref{3.26}. This completes the proof of Theorem \ref{th3.1}.

\end{document}